\documentclass[12pt, oneside, a4paper]{article}
\usepackage{graphicx}
\usepackage{tikz}
\usepackage{cite}
\usepackage{hyperref}
\usepackage{amsfonts}
\usepackage{fancyhdr,amsthm,amsmath,amssymb,mathrsfs}
\usepackage{bookmark}
\usepackage{enumerate}
\usetikzlibrary{intersections}

\newtheorem{thm}{\bf Theorem}[section]
\newtheorem{lem}[thm]{Lemma}

\newtheorem{cor}[thm]{Corollary}
\theoremstyle{definition}

\newtheorem{eg}{Example}
\numberwithin{equation}{section}
\newenvironment{prf}{\noindent {\bf Proof.}\rm}{\qed}

\newcommand{\DOI}[1]{DOI: \href{https://doi.org/#1}{#1}}

\begin{document}
\title{\Large\bf Some properties of congruence lattices of path semigroups\thanks{The first author is supported by the China Scholarship Council and the Postgraduate Research $\&$ Practice Innovation Program of Jiangsu Province (Grant No. KYCX22-1533). The third author is supported by the National Natural Science Foundation of China (Grant No. 11771212).}}
\author{Yongle Luo$^{a}$, Zhengpan Wang$^{b}$ and Jiaqun Wei$^{a}$\thanks{Corresponding author} \\
{\small\em $^{a}$School of Mathematical Sciences, Nanjing Normal University, $210023$ Nanjing, China} \\
{\small\em $^{b}$School of Mathematics and Statistics, Southwest University, $400715$ Chongqing, China} \\
{\small\em E-mail: yongleluo@nnu.edu.cn, \; zpwang@swu.edu.cn, \; weijiaqun@njnu.edu.cn} }
\date{}
\maketitle

\begin{abstract}
Each quiver corresponds to a path semigroup, and such a path semigroup also corresponds to an associative $K$-algebra over an algebraically closed field $K$. Let $Q$ be a quiver and $S_Q$, $KQ$ be its path semigroup, path algebra, respectively. In this paper, we study some properties of the congruence lattice of $S_Q$. First, we show that there is a one-to-one correspondence between congruences on $S_Q$ and certain algebraic ideals of $KQ$. Based on such a description, we consider acyclic quivers and show that the congruence latices of such path semigroups are strong upper semimodular but not necessarily lower semimodular. Moreover, we provide some equivalent conditions for the congruence lattices to be modular and distributive.
\medskip

{\bf 2020 Mathematics Subject Classification:} 20M10, 16D25, 06C99, 06D99, 05C20

{\bf Keywords}: path semigroups; path algebras; congruence lattices; algebraic ideals; modular lattices; distributive lattices.
\end{abstract}

\section{Introduction and main results}

Throughout this article, $K$ will denote an arbitrary algebraically closed field. Let $S$ be a semigroup. The {\it semigroup algebra} $K[S]$ is the $K$-algebra whose underlying $K$-vector space has a basis the set of all nonzero elements in $S$ and the multiplication extending the multiplication on $S$. Recall that an equivalence $\rho$ on $S$ is called a {\it congruence}, if $(x, y)\in \rho$ implies that $(ax, ay)\in \rho$ and $(xa, ya)\in \rho$ for all $x, y, a \in S$.  A congruence $\rho$ is called a {\it Rees congruence}, if $\rho$ is generated by a semigroup ideal $I$, that is, $\rho = (I\times I) \cup 1_S$. Congruences on semigroups are related to algebraic ideals of their semigroup algebras. Every congruence $\rho$ on $S$ determines an algebraic ideal $I_\rho$ of $K[S]$, where $I_\rho$ is the $K$-vector space with basis the set of all $x-y$ with $(x, y)\in \rho$. And the quotient algebra $K[S]/I_\rho$ is isomorphic to the semigroup algebra $K[S/\rho]$ , where $S/\rho$ is the quotient semigroup of $S$. And conversely, every algebraic ideal $I$ of $K[S]$ can generate a congruence $\rho_I$ on $S$, where $\rho_I =\{(x, y) \in S\times S \mid x-y \in I\}$.

Similar to the lattices of ideals of algebras, the poset of all congruences on $S$ partially ordered by containments is also a complete lattice, where the minimum is the identical relation $1_S$, and the maximum is the universal relation $S\times S$. Congruence theory is one of the main topic in semigroup theory. Congruences and congruence lattices of semigroups are well studied, one can refer to \cite{Mitsch,Mitsch2} as a general overview, or recent investigations, such as \cite{Anag,Luo}.

The path semigroup (see Subsection \ref{notions}) is a natural and important example. Let $Q$ be a quiver and $S_Q$ be its path semigroup. The semigroup algebra $K[S_Q]$ is just the path algebra which is well studied in representation theory of algebras. By convention we use $KQ$ rather than $K[S_Q]$ to denote the path algebra of $Q$. Algebraic ideals of path algebras are very important. For examples, a classical result of Gabriel shows that any finite dimensional $K$-algebra over an algebraically closed field $K$ is Morita equivalent to a quotient algebra $KQ/I$, where $I$ is an admissible ideal. More recently, deep relations between ideals of algebras and congruences of the lattice of torsion classes have been discovered by DIRRT \cite [Section 5.2] {Demonet}.

Path semigroups also play an important role in the study of graph inverse semigroups \cite{MesyanMitchell}, Leavitt path algebras \cite{AbramsPino}, partition theory \cite {Bergeleson} and partial semigroups \cite{East}. There are also many articles studied path semigroups by their own right. Zheng and Guo \cite{Zheng} used a categorical characterization of path semigroups, which showed that path semigroups are precisely the free objects in the category of locality semigroups with a rigid condition. Forsberg studied the effective dimensions of path semigroups over an uncountable field \cite{Forsberg}.
However, there are few known properties in congruence theory of path semigroups.

In the present article, we mainly studied some properties of congruence lattices of path semigroups. Let $Q$ be an arbitrary quiver and $S_Q$, $KQ$ be its path semigroup, path algebra, respectively. Denote by $\mathscr{C}(S_Q)$, $\mathcal{SI}(KQ)$ the set of all congruences on $S_Q$, the set of all {\it special algebraic ideals} of $KQ$ (see Subsection \ref{notions}), respectively. The two posets partially ordered by containments are both complete lattices. We proved that there is a closed relationship between the two lattices.

\begin{thm} [Theorem \ref{lattices iso}]
Let $Q$ be a quiver. Then the lattices $\mathscr{C}(S_Q)$ and $\mathcal{SI}(KQ)$ are isomorphic.
\end{thm}

Recall that an algebra $A$ is {\it hereditary}, if every submodule of a projective $A$-module is also projective. This kind of algebras has aroused great interests in representation theory. It turns out that a finite dimensional $K$-algebra $A$ over an algebraic closed field $K$ is basic, connected and hereditary if and only if $A\cong KQ$, where $Q$ is a finite, connected and acyclic quiver \cite[Theorem 1.7 of Chapter \uppercase\expandafter{\romannumeral7}]{Assem}. In the present article, we mainly considered the distributivity, modularity, and semimodularity (see Subsection \ref{lattice preliminary}) of $\mathscr{C}(S_Q)$ when $Q$ is acyclic.

\begin{thm}[Theorem \ref{congruences upper semimodular}]
Let $Q$ be an acyclic quiver. Then $\mathscr{C}(S_Q)$ is strong upper semimodular but not necessarily lower semimodular.
\end{thm}

Besides, we provided some necessary and sufficient conditions for $\mathscr{C}(S_Q)$ to be lower semimodular and modular, as follows:

\begin{thm}[Theorem \ref{modular condition}]
Let $Q$ be an acyclic quiver. Then the following statements are equivalent.
\begin{enumerate}[$(1)$]
\item There are at most two paths from one vertex to another.
\item $\mathscr{C}(S_Q)$ is modular.
\item $\mathscr{C}(S_Q)$ is strong lower semimodular.
\item $\mathscr{C}(S_Q)$ is lower semimodular.
\end{enumerate}
\end{thm}

In addition, we also described necessary and sufficient conditions for $\mathscr{C}(S_Q)$ to be distributive.

\begin{thm}[Theorem \ref{distributive condition}]
Let $Q$ be an acyclic quiver. Then the following statements are equivalent.
\begin{enumerate}[$(1)$]
\item There is at most one path from one vertex to another.
\item Each congruence $\rho \in \mathscr{C}(S_Q)$ is a Rees congruence.
\item $\mathscr{C}(S_Q)$ is distributive.
\end{enumerate}
\end{thm}

A classical theorem of Gabriel shows that a connected, finite dimensional hereditary algebra $A$ is {\it representation-finite} (the number of the isomorphism classes of indecomposable finite dimensional $A$-modules is finite) if and only if $A \cong KQ$, where $Q$ is a quiver whose underlying graph is one of the Dynkin graphs\cite[Theorem 5.10 of Chapter \uppercase\expandafter{\romannumeral7}]{Assem}. As an interesting application, we get the following corollary.

\begin{cor}[Corrollary \ref{application tree}]
Let $Q$ be a quiver. If the underlying graph of $Q$ is a tree, in particular, where it is a Dynkin graph, then $\mathscr{C}(S_Q)$ is distributive.
\end{cor}

\section{Preliminaries}

\subsection{Lattices}\label{lattice preliminary}
In this subsection, we recall some lattice theory from \cite{Gratzer}.

Let $(L, \leq)$ be a poset and $a, b \in L$. We say that $a$ {\it covers} $b$, and write $a\succ b$ or $b \prec a$, if $a > b$ and there is no element $c\in L$ such that $a > c > b$.

A poset $L$ is called a {\it lattice} if, for all $a, b \in L$, there exists a least upper bound $a\vee b\in L$, called the {\it join} of $a$ and $b$, and a greatest lower bound $a\wedge b \in L$, called the {\it meet} of $a$ and $b$.
Recall, that a lattice $L=(L,\leq, \vee, \wedge)$ is
\begin{enumerate}[(i)]
\item {\it distributive}, if $(a\vee b)\wedge c =(a\wedge c)\vee (b\wedge c)$ for all $a, b, c \in L$.
\item {\it modular}, if $a \leq c$ implies that $(a\vee b)\wedge c = a\vee (b\wedge c)$ for all $a, b, c\in L$.
\item {\it strong upper semimodular}, if $a\succ a\wedge b$ implies that $a \vee b \succ b$ for all $a, b \in L$.
\item {\it strong lower semimodular}, if $a\vee b \succ a$ implies that $b\succ a\wedge b$ for all $a, b \in L$.
\item {\it upper semimodular}, if $a, b \succ a\wedge b$ together imply that $a\vee b \succ a, b$ for all $a, b \in L$.
\item {\it lower semimodular}, if $a\vee b \succ a, b$ together imply that $a, b \succ a\wedge b$ for all $a, b \in L$.
\end{enumerate}

It is clear that each distributive lattice is modular, and each modular lattice is both strong upper semimodular and strong lower semimodular. Every strong upper (lower) semimodular lattice is upper (lower) semimodular. But conversely, modular lattices are not always distributive. An infilite lattice which is both strong upper semimodular and strong lower semimodular is not necessarily modular. Moreover, the upper (lower) semimodularity does not imply the strong upper (lower) semimodularity in general.

\begin{figure}[htbp]
\centering
\begin{tikzpicture}[scale=0.7]
\fill(0,2.5)circle(2pt);
\fill (1,1)circle(2pt);
\fill (2,2)circle(2pt);
\fill (2,3)circle(2pt);
\fill (1,4)circle(2pt);
\draw (1,1)--(2,2);
\draw (2,2)--(2,3);
\draw (1,4)--(2,3);
\draw (0,2.5)--(1,4);
\draw (0,2.5)--(1,1);
\node[below] at (1,1) {$\mathfrak{N}_5$};

\fill(5,1)circle(2pt);
\fill (5,2.5)circle(2pt);
\fill (5,4)circle(2pt);
\fill (4,2.5)circle(2pt);
\fill (6,2.5)circle(2pt);
\draw (5,1)--(5,2.5);
\draw (5,2.5)--(5,4);
\draw (4,2.5)--(5,4);
\draw (6,2.5)--(5,4);
\draw (5,1)--(4,2.5);
\draw (5,1)--(6,2.5);
\node[below] at (5,1) {$\mathfrak{M}_3$};
\end{tikzpicture}
\caption{\footnotesize The lattices $\mathfrak{N}_5$ and $\mathfrak{M}_3$.}\label{Figure 1}
\end{figure}
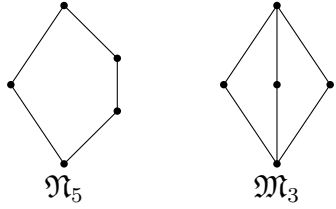

A sublattice $L_1$ of a lattice $L$ is called a {\it pentagon}, respectively a {\it diamond}, if $L_1$ is isomorphic to $\mathfrak{N}_5$, respectively to $\mathfrak{M}_3$, shown in Figure \ref{Figure 1}. The following two lemmas describle a lattice to be distributive or modular.

\begin{lem}\cite [Theorem 101]{Gratzer} \label{diagrams of distri}
A lattice $L$ is distributive if and only if it does not contain a pentagon $\mathfrak{N}_5$ or a diamond $\mathfrak{M}_3$.
\end{lem}

\begin{lem}\cite[Theorem 102]{Gratzer}\label{diagrams of modular}
A lattice $L$ is modular if and only if it does not contain a pentagon $\mathfrak{N}_5$.
\end{lem}

Let $L_1= (L_1, \leq_1)$ and $L_2= (L_2, \leq_2)$ be any two lattices. We say that $L_1$ and $L_2$ are {\it isomorphic} if, there is a bijection $\phi: L_1 \rightarrow L_2$ such that $a \leq_1 b$ in $L_1$ if and only if $a\phi \leq_2 b\phi$ in $L_2$.

\subsection{Special algebraic ideals of path algebras} \label{notions}

In this subsection, we give some concepts about special algebraic ideals of path algebras. For simplicity, we omit some notations about quivers and path algebras, which can be found in \cite{Assem}.

Let $Q =(Q_0, Q_1, s, t)$  be a quiver. The {\it path semigroup} $S_Q$ of $Q$ is the semigroup with zero generated by $Q_0$ and $Q_1$. More precisely, the multiplication in $S_Q$ is the same as the multiplication in the path algebra $KQ$. Nonzero elements in $S_Q$ are all paths in $Q$. And the set of all idempotents in $S_Q$ is $\{0\}\cup Q_0$.

Let $I$ be an algebraic ideal of $KQ$ and $A\subseteq I$, we say that $I$ is {\it generated} by $A$ if $I=<A>$, where $<A>$ is the smallest algebraical ideal containing $A$. Then, for any $x\in I$, $x$ is a finite sum of elements of the form $uav$, where $u, v \in KQ$ and $a\in A$. Since every element in $KQ$ is a finite sum of elements of the form $k\alpha$, where $k\in K$ and $\alpha$ is a path. Then $x$ can be expressed as follow:
$$ x= k_1u_1a_1v_1+...+k_mu_ma_mv_m,$$
where $k_i \in K$, $u_i, v_i$ are paths and $a_i \in A$ for any $1\leq i \leq m$.

Denote by $\mathcal{I}(KQ)$ the set of all algebraic ideals of $KQ$. A classical lattice result shows that $(\mathcal{I}(KQ), \subseteq, +, \cap)$ is a complete modular lattice. More precisely, let $\{I_j\}_{j\in J}$ be a family of algebraic ideals, where $J$ is an index set. The meet of $\{I_j\}_{j\in J}$ is $\bigcap\limits_{j\in J} I_j$, and the join is
\begin{equation}\label{equaltion 1}
\bigvee_{j\in J} I_j = \{x=\sum_{i=1}^{n} x_{j_i} \mid x_{j_i}\in I_{j_i},j_{i} \in J, n \in \mathbb{Z}^{+} \}
\end{equation}

For any path $\alpha$ in $Q$, we call the element $\rho = \alpha$ a {\it monomial relation} in $KQ$. For any paths $\alpha\neq \beta \in Q$ with $s(\alpha) = s(\beta)$ and $t(\alpha) = t(\beta)$, we call the element $\rho = \alpha - \beta$ a {\it commutative relation} in $KQ$ (Note that the concepts here are slightly different from those used in representation theory of algebras \cite{Assem}). In this case, we also define $s(\rho)= s(\alpha)$ and $t(\rho)= t(\alpha)$. Monomial relations and commutative relations are called {\it special relations}. Let $\rho$ be a special relation in $KQ$. Certainly, $\alpha\rho\beta=0$ or $\alpha\rho\beta$ is again a special relation for any paths $\alpha, \beta$.
An algebraic ideal is called {\it special}, if it is the zero ideal or it can be generated by special relations.

Denote by $\mathcal{SI}(KQ)$ the set of all special algebraic ideals of $KQ$. Note that the zero ideal and $KQ$ are both special algebraic ideals. Thus $\mathcal{SI}(KQ)$ is a bounded poset ordered by inclusions. Moreover, we can see that ($\mathcal{SI}(KQ), \subseteq)$ is also a complete lattice. Let $\{I_j\}_{j\in J}$ be a family of special algebraic ideals, where $J$ is an index set. It is obvious that $\bigvee\limits_{j\in J} I_j$ defined in formula (\ref{equaltion 1}) is special and thus it is also the join of $\{I_j\}_{j\in J}$ in $\mathcal{SI}(KQ)$. But in general, $\bigcap\limits_{j\in J} I_j$ is no longer generated by special relations (see Example \ref{counter not lower}). Hence $(\mathcal{SI}(KQ), \subseteq)$ is not a sublattice of $(\mathcal{I}(KQ), \subseteq, +, \cap)$. Indeed, the meet of $\{I_j\}_{j\in J}$ in $\mathcal{SI}(KQ)$ is given by
\begin{equation}\label{equation 2}
\bigwedge_{j\in J} I_j = <\rho \mid \rho ~\mbox{is a relation in}~\bigcap\limits_{j\in J} I_j>.
\end{equation}
That is, for any $x \in \bigwedge\limits_{j\in J} I_j$, $x$ can be expressed as $x= \sum\limits_{i=1}^{n} u_i \rho_i v_i$, where $n \in \mathbb{Z}^{+}$, each $u_i, v_i \in KQ$, and each $\rho_i$ is a relation in $\bigcap\limits_{j\in J} I_j$ for any $1 \leq i \leq n$.

Let $A$ be any subset of $KQ$. We denote by $span ~A$ the $K$-vector subspace generated by $A$.

\section{Proof of the main results}
First, we give some properties of congruences on path semigroups, which will be needed for the proof of our main results.

\begin{lem}\label{property of congruence}
Let $Q$ be a quiver. For any congruence $\rho$ in $\mathscr{C}(S_Q)$, we have the following.
\begin{enumerate}[$(1)$]
\item If $(\alpha, \beta) \in \rho$ and $\alpha \notin 0\rho$, then $s(\alpha) = s(\beta)$ and $t(\alpha) = t(\beta)$.
\item If $(e, f)\in \rho$, where $e, f$ are different vertices in $Q$, then $e, f\in 0\rho$.
\end{enumerate}
\end{lem}

\begin{prf} (1) If $(\alpha, \beta) \in \rho$, then it follows from $(s(\alpha)\alpha, s(\alpha)\beta)\in\rho$ that $(\alpha, s(\alpha)\beta) \in \rho$. Note that $s(\alpha)\beta \neq 0$ since $\alpha \notin 0\rho$. Hence, we conclude $s(\alpha)= s(\beta)$. By duality, we also have $t(\alpha)= t(\beta)$.

(2) Assume that $e \notin 0 \rho$. Then it follows from $(1)$ that $s(e) = s(f)$, that is, $e= f$, which is a contradiction.
\end{prf}\medskip

\begin{lem}\label{path beta and congruence}
Let $Q$ be a quiver. For any path $\beta$ and any congruence $\rho$ in $\mathscr{C}(S_Q)$, if there exist some elements $\alpha_1, ..., \alpha_{m+p}, \beta_{m+1}, ..., \beta_{m+p}$ in $S_Q$, where $m, p$ are both positive integers, such that the following hold:
\begin{enumerate}[$(1)$]
\item  $\beta =k_1\alpha_1 + ... +k_m \alpha_m + k_{m+1}(\alpha_{m+1}-\beta_{m+1})+...+k_{m+p}(\alpha_{m+p}-\beta_{m+p})$, where $k_i \in K\backslash\{0\}$ for any $1 \leq i \leq m+p$.
\item $\beta, \alpha_1, ..., \alpha_m$ are $K$-linearly independent.
\item $\alpha_{m+1}-\beta_{m+1}, ..., \alpha_{m+p}-\beta_{m+p}$ are $K$-linearly independent, and $(\alpha_{m+j},\beta_{m+j})\in \rho$ for any $1\leq j \leq p$.
\end{enumerate}
Then, there exists $x\in \{0, \alpha_1, ..., \alpha_m\}$ such that $(\beta, x)\in \rho $, and
\begin{equation}
\beta- x \in span~\{\alpha_{m+1}-\beta_{m+1}, ..., \alpha_{m+p}- \beta_{m+p}\}.\notag
\end{equation}
\end{lem}

\begin{prf} We proceed by induction on $p$. If $p=1$, that is,
\begin{equation}
\beta = k_1 \alpha_1 +... +k_m \alpha_m +k_{m+1}(\alpha_{m+1} -\beta_{m+1}),\notag
\end{equation}
then we obtain from condition $(2)$ that $\beta = \alpha_{m+1}$ or $\beta = \beta_{m+1}$.

If $\beta = \alpha_{m+1}$. Since $\alpha_{m+1} -\beta_{m+1}$ is $K$-linearly independent, then $\alpha_{m+1} \neq \beta_{m+1}$. So we have $k_{m+1}=1$ and $0= k_1 \alpha_1 +... +k_m \alpha_m -\beta_{m+1}$. Since $k_i \neq 0$ for all $1\leq i\leq m$, it follows from $(1)$ that $\beta_{m+1} \neq0$ and there exists $1\leq i \leq m$ such that $\beta_{m+1} =\alpha_i$. That is, $(\beta, \alpha_i)= (\alpha_{m+1}, \beta_{m+1})\in \rho$, and $\beta- \alpha_{i} \in span~ \{\alpha_{m+1}-\beta_{m+1}\}$. Similarly, if $\beta = \beta_{m+1}$, the conclusion also holds.

Next, we may assume that the conclusion holds for any $1 \leq p \leq n$. Let
\begin{equation} \label{EE1}
\beta = k_1 \alpha_1 +... +k_m \alpha_m +k_{m+1}(\alpha_{m+1} -\beta_{m+1}) + ... + k_{m+n+1}(\alpha_{m+n+1}-\beta_{m+n+1}).
\end{equation}

Note that $\alpha_{m+j} \neq \beta_{m+j}$ for any $1\leq j \leq n+1$. Without loss of generality, we may assume that $\beta= \alpha_{m+1}= ... = \alpha_{m+s}$ and $k_{m+1} +... +k_{m+s}=1$ for some $s$ with $1 \leq s \leq n+1$. By condition $(3)$, we get that $\beta_{m+1}, ..., \beta_{m+s}$ are different from each other.

{\bf Case 1}.  $\{0, \alpha_1, ..., \alpha_m\}\cap \{\beta_{m+1}, ..., \beta_{m+s}\} \neq \emptyset$.

Clearly, there exists $x \in \{0, \alpha_1, ..., \alpha_m\}$  and $\beta_{m+j}$ for some $j$ with $1\leq j \leq s$ such that $x=\beta_{m+j}$. Then we have $(\beta, x)= (\alpha_{m+j}, \beta_{m+j}) \in \rho$, and \begin{equation}
\beta- x\in span~\{\alpha_{m+1} - \beta_{m+1}, ..., \alpha_{m+n+1}- \beta_{m+n+1}\}.\notag
\end{equation}

{\bf Case 2}.  $\{0, \alpha_1, ..., \alpha_m\}\cap\{\beta_{m+1}, ..., \beta_{m+s}\}=\emptyset$.

In this case $\beta, \alpha_1, ..., \alpha_m, \beta_{m+1}, ..., \beta_{m+s}$ are different paths, that is, they are $K$-linearly independent. Then immediately we have $1\leq s\leq n$. Otherwise, it follows $s= n+1$ and formula (\ref{EE1}) that
\begin{equation}
0= k_1\alpha_1 +...+k_m\alpha_m-k_{m+1}\beta_{m+1}-...-k_{m+n+1}\beta_{m+n+1}.\notag
\end{equation}
Then $\alpha_1$, ..., $\alpha_m$, $\beta_{m+1}$, ..., $\beta_{m+n+1}$ are $K$-linearly dependent, which is a contradiction.

Note that
$$\begin{array}{rcl}
0 &= &k_1\alpha_1 +...+k_m\alpha_m - k_{m+1}\beta_{m+1}- ...-k_{m+s}\beta_{m+s}  \\
&& + k_{m+s+1}(\alpha_{m+s+1}-\beta_{m+s+1})+...+k_{m+n+1}(\alpha_{m+n+1}-\beta_{m+n+1}).
\end{array}$$
That is, $$\begin{array}{rcl}
\beta_{m+s}&=&\frac{1}{k_{m+s}}[k_1\alpha_1 +...+k_m\alpha_m - k_{m+1}\beta_{m+1} -... - k_{m+s-1}\beta_{m+s-1}  \\
&&~~~~~~+k_{m+s+1}(\alpha_{m+s+1}-\beta_{m+s+1})+...+k_{m+n+1}(\alpha_{m+n+1}-\beta_{m+n+1})].
\end{array}$$

By induction, there exists $x\in\{0, \alpha_1, ..., \alpha_m, \beta_{m+1}, ..., \beta_{m+s-1}\}$ such that $(\beta_{m+s}, x)\in \rho $, and
\begin{equation}\label{eq4}
\beta_{m+s}- x \in span~\{\alpha_{m+s+1}-\beta_{m+s+1}, ..., \alpha_{m+n+1}-\beta_{m+n+1}\}.
\end{equation}

Since $\beta= \alpha_{m+s}$ and $(\alpha_{m+s}, \beta_{m+s})\in \rho$, we have $(\beta, x)\in \rho$. Moreover,
\begin{equation}
\beta - x =(\beta - \beta_{m+s}) + (\beta_{m+s}-x)= (\alpha_{m+s} - \beta_{m+s})+(\beta_{m+s}-x).\notag
\end{equation}
and thus $\beta - x \in span~\{\alpha_{m+1}-\beta_{m+1}, ..., \alpha_{m+n+1}-\beta_{m+n+1}\}$.

It remains to prove $x \in \{0, \alpha_1, ..., \alpha_m\}$. Indeed, if $x=\beta_{m+j}$ for some $1\leq j \leq s-1$, then it follows from $\alpha_{m+s} = \alpha_{m+j}$,
\begin{equation}
-(\alpha_{m+s}- \beta_{m+s})+(\alpha_{m+j} -\beta_{m+j})= \beta_{m+s}-\beta_{m+j} = \beta_{m+s}- x,\notag
\end{equation}
and formula (\ref{eq4}) that $\alpha_{m+1}-\beta_{m+1}, ..., \alpha_{m+n+1}-\beta_{m+n+1}$ are $K$-linearly dependent, which is a contradiction. So we conclude $x \in \{0, \alpha_1, ..., \alpha_m\}$. This finishes the proof.
\end{prf}\medskip

\begin{lem}\label{path beta in 0rho}
Let $Q$ be a quiver. For any path $\beta$ and any congruence $\rho$ in $\mathscr{C}(S_Q)$, if there exist some elements $\alpha_1, ..., \alpha_{m}, \beta_{1}, ..., \beta_{m}$ in $S_Q$, where $m$ is a positive integer , such that the following hold:

\begin{enumerate}[$(1)$]
\item $\beta = k_1 (\alpha_1 -\beta_1) + ... +k_m(\alpha_{m} -\beta_m)$, where $k_i \in K\backslash\{0\}$ for any $1 \leq i \leq m$.
\item $\alpha_1 -\beta_1, ..., \alpha_m -\beta_m$ are $K$-linearly independent, and $(\alpha_{i},\beta_{i})\in \rho$ for any $1\leq i \leq m$.
\end{enumerate}
Then $(\beta, 0)\in \rho$.
\end{lem}

\begin{prf} We proceed by induction on $m$. If $m=1$, that is, $\beta = k_1(\alpha_1 - \beta_1)$, then we easily have $\beta = \alpha_1$ or $\beta=\beta_1$. For the former case, $k_1 =1$ and $\beta_1=0$. For the latter case, $k_1 =-1$ and $\alpha_1 =0$. Certainly, the conclusion holds.

Suppose that the conclusion holds for any $1\leq m \leq n$. Now we consider the case when $\beta = k_1 (\alpha_1 -\beta_1) + ... +k_{n+1}(\alpha_{n+1} -\beta_{n+1})$. Since $\alpha_i \neq \beta_i$ for any $1\leq i \leq n+1$. We may assume that $\beta= \alpha_{1}= ... = \alpha_{s}$ and $k_{1} +... +k_{s} = 1$ for some $s$ with $1 \leq s \leq n+1$. Thus, $\beta_1, ..., \beta_s$ are different elements from each other. If $0\in \{\beta_1, ..., \beta_s\}$, there is nothing to prove. So we assume that $\beta_1, ..., \beta_s$ are different paths. Note that
$$\begin{array}{rcl}
0=-k_1\beta_1 -...-k_s\beta_s +k_{s+1}(\alpha_{s+1}-\beta_{s+1})+...+k_{n+1}(\alpha_{n+1} - \beta_{n+1}),
\end{array}$$
that is,
$$\begin{array}{rcl}
\beta_s = \frac{1}{k_s}[-k_1\beta_1 - ... - k_{s-1}\beta_{s-1} + k_{s+1}(\alpha_{s+1}-\beta_{s+1})+...+k_{n+1}(\alpha_{n+1} - \beta_{n+1})].
\end{array}$$
Clearly, we have $s\leq n$. Indeed, if $s=n+1$, then we get from $$0=-k_1\beta_1-... -k_{n+1}\beta_{n+1}$$ and $k_i \neq0$ for all $1\leq i \leq n+1$ that $\beta_1, ..., \beta_s$ are $K$-linearly dependent, which is a contradiction.

If $s=1$, then $(\beta_1, 0)\in \rho$ by induction. Note that $\beta = \alpha_1$ and $(\alpha_1, \beta_1) \in \rho$, then we have $(\beta, 0)\in \rho$.

It remains to consider the case when $2\leq s \leq n$. By Lemma \ref{path beta and congruence}, there exists $x\in \{0, \beta_1, ..., \beta_{s-1}\}$ such that $(\beta_s, x)\in \rho$ and $$\beta_s - x \in span ~\{\alpha_{s+1}-\beta_{s+1}, ..., \alpha_{n+1}-\beta_{n+1}\}.$$
We claim that $x=0$. Indeed, if $x=\beta_j$ for some $j$ with $1\leq j \leq s-1$, since $\alpha_s = \alpha_j $, then it follows
$$-(\alpha_{s}- \beta_{s})+(\alpha_j -\beta_j)= \beta_{s}-\beta_j = \beta_{s}- x.$$
Note that $\beta_{s}- x\in span~\{\alpha_{s+1}-\beta_{s+1}, ..., \alpha_{n+1}-\beta_{n+1}\}$. Thus $\alpha_{1}-\beta_{1}, ..., \alpha_{n+1}-\beta_{n+1}$ are $K$-linearly dependent, which is a contradiction. So we conclude that $x=0$. Consequently, it follows from $\beta = \alpha_s$, $(\beta_s, 0)\in \rho$ and $(\alpha_s, \beta_s)\in \rho$ that $(\beta, 0)\in \rho$.
\end{prf}\medskip

\begin{lem}\label{paths alpha and beta in rho}
Let $Q$ be a quiver. For any elements $\alpha, \beta$ in $S_Q$ and any congruence $\rho$ in $\mathscr{C}(S_Q)$, if there exist some elements $\alpha_1, ..., \alpha_{m}, \beta_{1}, ..., \beta_{m}$ in $S_Q$, where $m$ is a positive integer and $(\alpha_i, \beta_i)\in \rho$ for all $1 \leq i \leq m$, such that
$$\alpha - \beta = k_1 (\alpha_1 -\beta_1) + ... +k_m(\alpha_{m} -\beta_m),$$
where $k_i \in K$ for all $1 \leq i \leq m$. Then $(\alpha, \beta)\in \rho$.
\end{lem}

\begin{prf} If $\alpha = \beta$, there is noting to show. If $\alpha = 0$ and $\beta \neq 0$, Without loss of generality, we may assume that $k_i\neq 0$ for any $1\leq i\leq m$ and $\alpha_1 - \beta_1, ..., \alpha_m- \beta_m$ are $K$-linearly independent. It follows from Lemma \ref{path beta in 0rho} that $(\beta, 0)\in \rho$. Thus, $(\alpha, \beta)\in \rho$. If $\beta =0$ and $\alpha \neq 0$, the proof is similar.

Now, we consider the case when $\alpha$ and $\beta$ are different paths.
We can also assume that $k_i \neq 0$ for any $1\leq i\leq m$, $\alpha_1 -\beta_1$, ..., $\alpha_m -\beta_m$ are $K$-linearly independent.
Note that
$$\alpha = \beta + k_1 (\alpha_1 -\beta_1) + ... +k_m(\alpha_{m} -\beta_m).$$
By Lema \ref{path beta and congruence}, then we have $(\alpha, \beta)\in \rho$ or $(\alpha, 0)\in \rho$. If $(\alpha, 0)\in \rho$, then $\alpha \in span~\{\alpha_1 -\beta_1, ..., \alpha_m -\beta_m\}$. In this case, $\beta = \alpha -(\alpha-\beta) \in span~\{\alpha_1 -\beta_1, ..., \alpha_m -\beta_m\}$. By Lemma \ref{path beta in 0rho} again, we also have $(\beta, 0)\in \rho$. Hence we always have $(\alpha, \beta)\in\rho $.
\end{prf}\medskip

Next, we study the relationship between the lattice $\mathscr{C}(S_Q)$ and the lattice $\mathcal{SI}(KQ)$.

\begin{thm} \label{lattices iso}
Let $Q$ be a quiver. Then the lattices $\mathscr{C}(S_Q)$ and $\mathcal{SI}(KQ)$ are isomorphic.
\end{thm}
\begin{prf}
For any congruence $\rho \in \mathscr{C}(S_Q)$, $\rho$ determines an algebraic ideal
\begin{equation}
I_\rho = <\alpha-\beta \mid (\alpha, \beta)\in \rho>\notag
\end{equation}of $KQ$. It is easy to check that $I_\rho$ is a special algebraic ideal. To see this, it is enough to show that each generator defined in $I_\rho$ can be generated by special relations. Indeed, for any $(\alpha, \beta) \in \rho$, if $(\alpha ,0)\in \rho$, then $(\beta, 0) \in \rho$. Thus the monomials $\alpha = \alpha -0$ and $\beta = \beta - 0$ are generators of $I_\rho$, which together imply that $\alpha - \beta$ is generated by special relations. If $(\alpha, 0) \notin \rho$, then it follows from Lemma \ref{property of congruence} (1) that $s(\alpha)= s(\beta)$, $t(\alpha) = t(\beta)$. Thus $\alpha- \beta$ is a commutative relation in $I_\rho$.

Conversely, for any special algebraic ideal $I\in \mathcal{SI}(KQ)$, $I$ also determines a congruence $\rho_I \in \mathscr{C}(S_Q)$, where
$$\rho_I = \{(\alpha, \beta) \in S_Q \times S_Q\mid \alpha - \beta \in I\}.$$
Indeed, for any $\alpha \in S_Q$, since $\alpha- \alpha = 0 \in I$, then $(\alpha, \alpha)\in \rho_I$. For any $(\alpha, \beta) \in \rho_I$, by the construction of $\rho_I$ we have $\alpha- \beta \in I$. Then it follows $\beta - \alpha \in I$ and $(\beta, \alpha)\in \rho_I$. For any $(\alpha, \beta)\in \rho_I$ ans $(\beta, \gamma) \in \rho_{I}$, we get from $\alpha- \beta \in I$ and $\beta- \gamma \in I$ that $\alpha - \gamma = (\alpha- \beta) + (\beta - \gamma) \in I$. Hence we have $(\alpha, \gamma)\in \rho_I$. Next, we need to prove that $\rho_I$ satisfies the left and right compatible conditions. For any $(\alpha, \beta)\in \rho_I$ and any $\gamma \in S_Q$, it is easy to get from $\alpha- \beta \in I$ that $\gamma\alpha - \gamma\beta = \gamma (\alpha- \beta)\in I$, $\alpha\gamma - \beta\gamma = (\alpha- \beta)\gamma\in I$. Hence we have $(\gamma\alpha, \gamma\beta)\in \rho_I$, and  $(\alpha\gamma, \beta\gamma)\in \rho_I$.

Next, we prove that there is a one-to-one correspondence between $\mathscr{C}(S_Q)$ and $\mathcal{SI}(KQ)$. By the construction above, it is enough to show that for any $\rho \in \mathscr{C}(S_Q)$, $\rho= \rho_{I_{\rho}}$, and for any $I\in \mathcal{SI}(KQ)$, $I= I_{\rho_{I}}$.

Let $\rho$ be any congruence and $(\alpha, \beta)\in \rho$. Then we clearly have $\alpha- \beta \in I_\rho$. By the definition of $\rho_{I_{\rho}}$, certainly $(\alpha, \beta) \in \rho_{I_{\rho}}$. Thus, we proved $\rho \subseteq \rho_{I_{\rho}}$. Conversely, for any $(\alpha, \beta) \in \rho_{I_{\rho}}$, by the definition, we have $\alpha- \beta \in I_\rho$. Then $\alpha-\beta$ can be expressed as the following
$$\alpha- \beta= k_1u_1(\alpha_1 - \beta_1)v_1+...+k_nu_n(\alpha_n- \beta_n)v_n,$$
where $n$ is a positive integer, $k_i \in K$, $u_i, v_i$ are paths in $Q$, and $(\alpha_i, \beta_i)\in \rho $ for any $1\leq i \leq n$. By the transitivity of $\rho$, we have $(u_i\alpha_iv_i, u_i\beta_iv_i)\in \rho$. It follows from Lemma \ref{paths alpha and beta in rho} that $(\alpha, \beta)\in \rho$. Thus, we proved $\rho_{I_\rho} \subseteq \rho$.

Let $I$ be any special algebraic ideal. Since $I$ is generated by special relations, to prove $I \subseteq I_{\rho_{I}}$, we need only to prove that every special relation in $I$ is contained in $I_{\rho_I}$. For any  monomial relation $x \in I$, since $x= x-0\in I$, we have $(x, 0) \in \rho_I$, thus we get $x \in I_{\rho_I}$. For any commutative relation $x-y \in I$, we have $(x, y)\in \rho_I$, thus we have $x- y \in I_{\rho_I}$. Hence, we proved $I \subseteq I_{\rho_I}$. Conversely, for any generator $\alpha-\beta\in I_{\rho_I}$, by the definition, we have $(\alpha, \beta)\in \rho_I$. Thus, $\alpha- \beta \in I$. So we proved $I_{\rho_I}\subseteq I$.

It is obvious that the two constructions above both preserve partial orders. Thus the two lattices are isomorphic.
\end{prf}\medskip

\begin{eg}\label{eg1}
Let $Q$ be the quiver in Figure \ref{Figure 2}.

\begin{figure}[htbp]
\centering
\begin{tikzpicture} [scale=0.7]
\draw (1,1) circle(0.06);
\draw (3,1) circle(0.06);
\draw [-stealth] (1.1,1)--(2.9,1);
\node[below] at (1,1) {$1$};
\node[below] at (3,1) {$2$};
\node[above] at (2,1) {$\alpha$};
\end{tikzpicture}
\caption{\footnotesize An example.}\label{Figure 2}
\end{figure}

The path semigroup $S_Q=\{0, e_1, e_2, \alpha\}$. There are $5$ congruences on $S_Q$: $\rho_1 = 1_{S_Q}$, $\rho_2$ is the Rees congruence generated by $\{0, \alpha\}$, $\rho_3$ is the Rees congruence generated by $\{0, e_1, \alpha\}$, $\rho_4$ is the Rees congruence generated by $\{0, e_2,\alpha\}$, and $\rho_5$ is the universal congruence $S_Q\times S_Q$.

The corresponding special algebraic ideals in $\mathcal{SI}(KQ)$ are $I_1 = \{0\}$, $I_2 =<\alpha>$, $I_3 =<e_1>=span~\{e_1, \alpha\}$, $I_4 = <e_2>=span~\{e_2, \alpha\}$, and $I_5= <e_1, e_2>=KQ$.

The Hasse diagram of the lattice $\mathscr{C}(S_Q)$ ($\mathcal{SI}(KQ)$) is shown in Figure \ref{Figure 3}.

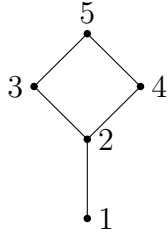
\begin{figure}[htbp]
\centering
\begin{tikzpicture}[scale=0.7]
\fill (1,1.5)circle(2pt); \node[right] at (1,1.5) {$1$};
\fill (1,3)circle(2pt); \node[right] at (1,3) {$2$};
\fill (1,5)circle(2pt); \node[above] at (1,5) {$5$};
\fill (0,4) circle(2pt); \node[left] at (0,4) {$3$};
\fill (2,4) circle(2pt);   \node[right] at (2,4) {$4$};
\draw (1,1.5)--(1,3);
\draw (1,3)--(2,4);
\draw (1,3)--(0,4);
\draw (0,4)--(1,5);
\draw (2,4)--(1,5);
\end{tikzpicture}
\caption{\footnotesize The Hasse diagram of the lattice $\mathscr{C}(S_Q)$ ($\mathcal{SI}(KQ)$)in Example \ref{eg1}.}\label{Figure 3}
\end{figure}

\end{eg}

In the representation theory of associative $K$-algebra, hereditary algebras play an important role. It is well known that a basic, connected, finite dimensional $K$-algebra $A$ is hereditary if and only if $A\cong KQ$, where $Q$ is a finite, connected, and acyclic quiver \cite[Theorem 1.7 of Chapter \uppercase\expandafter{\romannumeral7}]{Assem}. Next, we study the congruence lattices of acyclic quivers (not necessarily finite).

\begin{lem} \label{cover lemma 1}
Let $Q$ be an acyclic quiver. Let $I_1$ and $I_2$ be any special algebraic ideals of $KQ$ with $I_1 \subsetneq  I_2$. For any special relation $x\in I_2 \backslash I_1$, if there exist $\alpha_0, \beta_0 \in S_Q$ such that $\alpha_0 x \beta_0 \neq x $ and $ \alpha_0 x \beta_0 \notin I_1$, then $I_1 \subsetneq I_1 + <\alpha_0 x \beta_0> \subsetneq I_2$.
\end{lem}
\begin{prf} Note that $I_1 \subsetneq I_1 + <\alpha_0 x \beta_0>$ since $\alpha_0 x \beta_0 \notin I_1$.  To prove $I_1 + <\alpha_0 x \beta_0> \subsetneq I_2$, it is sufficient to prove $x\notin I_1 + <\alpha_0 x \beta_0>$. Suppose $x\in I_1 + <\alpha_0 x \beta_0>$. Without loss of generality, we may assume that
\begin{equation}
x = k_1\alpha_1 + ... + k_m \alpha_m + k_{m+1}\beta_1 +... +k_{m+p} \beta_p,\notag
\end{equation}
where $m, p$ are non-negative integers such that $m+p \geq 1$, $k_i \in K\backslash\{0\}$ for all $1\leq i \leq m+p$, $\alpha_1, ..., \alpha_m$ are special relations in $I_1$, and $\beta_1, ... , \beta_p$ are special relations in $<\alpha_0 x \beta_0>$. We know from $x\notin I_1$ that $p \geq 1$.
Moreover, we have
$$ x=s(x)xt(x)= s(x)(k_1\alpha_1 + ... + k_m \alpha_m)t(x) + s(x)(k_{m+1}\beta_1 +... +k_{m+p} \beta_p) t(x).$$
Since $s(x)(k_1\alpha_1 + ... + k_m \alpha_m)t(x) \in I_1$, we have $s(x)(k_{m+1}\beta_1 +... +k_{m+p} \beta_p) t(x)\neq 0$. Then there exists $\beta_i$ for some $1\leq i\leq p$ such that $s(x)\beta_{m+i}t(x)\neq 0$. Note that $\beta_{m+i}$ is a special relation generated by $\alpha_0 x \beta_0$. We may assume that $\beta_{m+i} = u \alpha_0 x \beta_0 v$, where $u, v$ are paths. Now we deduce $s(x)= s(\beta_{m+i})= s(u \alpha_0)= r(u\alpha_0)$. Since $Q$ is acyclic, we obtain that $u= \alpha_0 = s(x)$, similarly, $v=\beta_0 = r(x)$. Hence, we get a contradiction that $x= \alpha_0 x \beta_0$.
\end{prf}\medskip

\begin{lem} \label{covering}
Let $Q$ be an acyclic quiver. For any special algebraic ideals $I_1, I_2$ of $KQ$, if $I_1 \prec I_2$, then $I_2 = I_1 +span~\{x\}$ for any special relation $x\in I_2 \backslash I_1$.
\end{lem}

\begin{prf} Let $x$ be any relation in $I_2 \backslash I_1$. It is clear that $I_1 \subsetneq I_1 +span~\{x\} \subseteq I_2$ as $K$-vector spaces. Next, we show that $I_1 + span~\{x\}$ is a special algebraic ideal of $KQ$. Indeed, if there exist $\alpha_0, \beta_0 \in S$ such that $\alpha_0 x \beta_0\neq x$ and $\alpha_0 x \beta_0 \notin I_1$, by Lemma \ref{cover lemma 1}, then we know that $I_1 \subsetneq I_1 +<\alpha_0 x \beta_0> \subsetneq I_2$, which is a contradiction. Hence, we proved that $I_1 + span~\{x\}$ is an ideal for any special relation $x\in I_2 \backslash I_1$. Note that $I_1 \prec I_2$ and $I_1 +span~\{x\}$ is also generated by special relations. Certainly we have $I_2 = I_1 +span~\{x\}$.
\end{prf}\medskip

Next, we study the semimodularity of the congruence lattices of path semigroups of acyclic quivers.

\begin{lem}\label{special ideals upper semimodular}
Let $Q$ be an acyclic quiver. Then the lattice $(\mathcal{SI}(KQ), \subseteq, \vee, \wedge)$ is strong upper semimodular.
\end{lem}

\begin{prf}
Let $I_1$, $I_2$ be any two special ideals of $KQ$ with $I_1 \succ I_1\wedge I_2$. Then it follows $I_1 \nsubseteq I_2$. If $I_2 \subsetneq I_1$, that is, $I_1 \wedge I_2 = I_2$, then we clearly have $I_1 \vee I_2 = I_1 \succ I_2$. Next, we consider the case when $I_1$ and $I_2$ are not comparable.

Let $x$ be any special relation in $I_1 \backslash I_2$. It follows from Lemma \ref{covering} that $I_1= I_1 \wedge I_2 + span~\{x\}$. Let $Y$ be the set of all special relations in $I_2 \backslash I_1$. Then we can see that $I_2 = I_1 \wedge I_2 + span ~Y$. Indeed, since $Y \subseteq I_2$, we have $span~ Y \subseteq I_2$, and then $I_1\wedge I_2 + span~ Y \subseteq I_2$. On the other hand, for any special relation $\alpha \in I_2$, if $\alpha \in I_1 \cap I_2$, note that $\alpha$ is a special relation, then $\alpha \in I_1 \wedge I_2$. If $\alpha \notin I_1$, then $\alpha \in Y$. So we proved that any special relation $\alpha\in I_2$, $\alpha \in I_1\wedge I_2 + span ~Y$. Hence, $I_2 \subseteq I_1 \wedge I_2 + span~ Y$.

By the definition of the join of $I_1$ and $I_2$, we get $$I_1 \vee I_2 = I_1+I_2 = I_1 \wedge I_2 + span~\{x\}+ span~Y.$$
It is clear that $I_1 \vee I_2 \succ I_2$.
Hence we proved that the lattice $(\mathcal{SI}(KQ), \subseteq, \vee, \wedge)$ is strong upper semimodular.
\end{prf}\medskip

However, the lattice $(\mathcal{SI}(KQ), \subseteq, \vee, \wedge)$ is not necessarily lower semimodular as shown in the following example.

\begin{eg}\label{counter not lower}
Let $Q$ be the quiver in Figure \ref{Figure 4}.

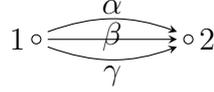
\begin{figure}[htbp]
\centering
\begin{tikzpicture}
\draw (1,0) circle(0.06);
\draw (3,0) circle(0.06);
\node[left] at (1,0) {$1$};  \draw[-stealth] (1.15, 0.1) to [bend left = 20] (2.85,0.1);\draw[-stealth] (1.15, -0.1) to [bend right = 20] (2.85,-0.1);
\node[right] at (3,0) {$2$};  \draw [-stealth] (1.15,0)--(2.85,0);
\node[above] at (2,0.2) {$\alpha$}; \node[above] at (2,-0.3) {$\beta$}; \node[below] at (2,-0.2) {$\gamma$};
\end{tikzpicture}
\caption{\footnotesize An example where $\mathcal{SI}(KQ)$ is not lower semimodular.}\label{Figure 4}
\end{figure}

The path semigroup $S_Q= \{0, e_1, e_2, \alpha, \beta, \gamma\}$. There are 18 special algebraic ideals of $KQ$: $I_1=\{0\}$, $I_2= <\alpha>$, $I_3= <\beta>$, $I_4= <\gamma>$, $I_5= <\alpha, \beta>$, $I_6= <\beta, \gamma>$, $I_7= <\alpha, \gamma>$, $I_8= <\alpha, \beta, \gamma>$, $I_9= <\alpha- \beta>$, $I_{10}= <\beta- \gamma>$, $I_{11}= <\alpha- \gamma>$,  $I_{12}= <\alpha, \beta- \gamma>$, $I_{13}= <\beta, \alpha- \gamma>$, $I_{14}= <\gamma, \alpha-\beta>$, $I_{15}= <\alpha- \beta,\alpha-\gamma>=<\alpha-\beta, \beta- \gamma>=<\alpha- \gamma, \beta- \gamma>$, $I_{16}= <e_1>=span~\{e_1, \alpha, \beta, \gamma\}$, $I_{17}= <e_2>=span~\{e_2, \alpha, \beta, \gamma\}$, $I_{18}= <e_1, e_2>=KQ$.

The Hasse diagram of the lattice $(\mathcal{SI}(KQ), \subseteq, \vee, \wedge)$ is shown in Figure \ref{Figure 5}. It is strong upper semimodular but not lower semimodular. Indeed, $I_{12}= span~\{\alpha, \beta- \gamma\}$, $I_{14}= span~\{\gamma, \alpha-\beta\}$, we have $I_{12} \vee I_{14}= span~\{\alpha, \beta, \gamma\}= I_8$.  Note that $I_{12} \cap I_{14}=span~\{\alpha-\beta+\gamma\}$ has no relation, so we have $I_{12}\wedge I_{14}=\{0\}$. It is clearly that $I_8 \succ I_{12}, I_{14}$, but there exists a special algebraic ideal $I_2= span~\{\alpha\}$ such that $\{0\} \subsetneq I_2 \subsetneq I_{12}$. Thus, the lattice is not lower semimodular.

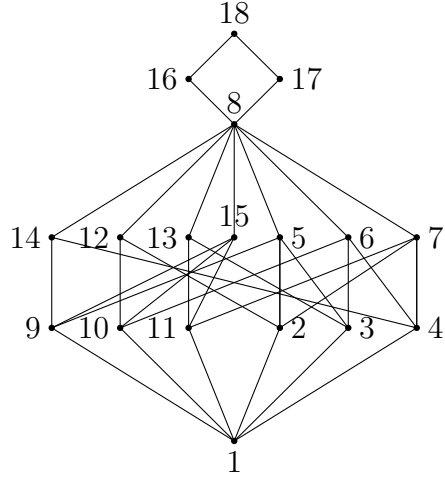
\begin{figure}[htbp]
\centering
\begin{tikzpicture} [scale=0.6]
\fill (0,0)circle(2pt); \node[below] at (0,0) {$1$};
\fill (-4,2.5)circle(2pt); \node[left] at (-4,2.5) {$9$};
\fill (-2.5,2.5)circle(2pt); \node[left] at (-2.5,2.5) {$10$};
\fill (-1,2.5)circle(2pt); \node[left] at (-1,2.5) {$11$};
\fill (1,2.5)circle(2pt); \node[right] at (1,2.5) {$2$};
\fill (2.5,2.5)circle(2pt); \node[right] at (2.5,2.5) {$3$};
\fill (4,2.5)circle(2pt); \node[right] at (4,2.5) {$4$};

\fill (-4,4.5)circle(2pt);  \node[left] at (-4,4.5) {$14$};
\fill (-2.5,4.5)circle(2pt);  \node[left] at (-2.5,4.5) {$12$};
\fill (-1,4.5)circle(2pt);  \node[left] at (-1,4.5) {$13$};
 \fill (0,4.5)circle(2pt); \node[above] at (0,4.5) {$15$};
\fill (1,4.5)circle(2pt); \node[right] at (1,4.5) {$5$};

\fill (0,7)circle(2pt); \node[above] at (0,7) {$8$};
\fill (2.5,4.5)circle(2pt);   \node[right] at (2.5,4.5) {$6$};
\fill (4,4.5)circle(2pt);\node[right] at (4,4.5) {$7$};

\fill (-1,8)circle(2pt); \node[left] at (-1,8) {$16$};
\fill (1,8)circle(2pt); \node[right] at (1,8) {$17$};
\fill (0,9)circle(2pt); \node[above] at (0,9) {$18$};

\draw (0,0)--(-4,2.5); \draw (0,0)--(-2.5,2.5); \draw (0,0)--(-1,2.5); \draw (0,0)--(1,2.5); \draw (0,0)--(2.5,2.5); \draw (0,0)--(4,2.5);
\draw (-4,4.5)--(-4,2.5); \draw (-2.5,4.5)--(-2.5,2.5); \draw (-1,4.5)--(-1,2.5); \draw (1,2.5)--(1,4.5); \draw (4,2.5)--(4,4.5);

\draw (-4,2.5)--(0,4.5); \draw (-4,2.5)--(1,4.5); \draw (-2.5,2.5)--(0,4.5); \draw (-2.5,2.5)--(2.5,4.5);  \draw (-1,2.5)--(0,4.5);  \draw (-1,2.5)--(4,4.5);
\draw (1,2.5)--(-2.5,4.5); \draw (1,2.5)--(4,4.5); \draw (1,2.5)--(1,4.5);  \draw (2.5,2.5)--(2.5,4.5); \draw (2.5,2.5)--(-1,4.5); \draw (2.5,2.5)--(1,4.5);
\draw (4,2.5)--(2.5,4.5); \draw (4,2.5)--(4,4.5); \draw (4,2.5)--(-4,4.5);
\draw (0,7)--(1,4.5); \draw (0,7)--(2.5,4.5); \draw (0,7)--(4,4.5);\draw (0,7)--(-2.5,4.5); \draw (0,7)--(-1,4.5); \draw (0,7)--(0,4.5);
\draw (0,7)--(-4,4.5); \draw (0,7)--(-1,8); \draw (0,7)--(1,8); \draw (1,8)--(0,9);\draw (-1,8)--(0,9);

\end{tikzpicture}
\caption{\footnotesize The Hasse diagram of the lattice $\mathcal{SI}(KQ)$ in Example \ref{counter not lower}.}\label{Figure 5}
\end{figure}
\end{eg}

By Theorem \ref{lattices iso}, Lemma \ref{special ideals upper semimodular} and Example \ref{counter not lower},  we immediately obtain the following theorem.

\begin{thm}\label{congruences upper semimodular}
Let $Q$ be an acyclic quiver. Then $\mathscr{C}(S_Q)$ is strong upper semimodular but not necessarily lower semimodular.
\end{thm}

The following theorem gives some equivalent characterizations of the congruence lattice to be lower semimodular.

\begin{thm}\label{modular condition}
 Let $Q$ be an acyclic quiver. Then the following statements are equivalent.
\begin{enumerate}[$(1)$]
\item There are at most two paths from one vertex to another.
\item $\mathscr{C}(S_Q)$ is modular.
\item $\mathscr{C}(S_Q)$ is strong lower semimodular.
\item $\mathscr{C}(S_Q)$ is lower semimodular.
\end{enumerate}
\end{thm}

\begin{prf} $(1)\Rightarrow (2)$. By Theorem \ref{lattices iso}, we should prove $\mathcal{SI}(KQ)$ is modular. Let $I_1$, $I_2$ be any special algebraic ideals of $KQ$. Next, we will show $I_1 \wedge I_2 = I_1 \cap I_2$, that is, $\mathcal{SI}(KQ)$ is a sublattice of $\mathcal{I}(KQ)$. In other words, every element $x$ in $ I_1 \cap I_2$, $x$ can be generated by special relations.

Note that there are at most two paths between two vertices. For any commutative relations $\alpha - \beta$ and $\alpha - \gamma$, clearly, we have $\beta = \gamma$. Moreover, we also should note that $<\alpha, \alpha-\beta> = <\alpha, \beta>$. Without loss of generality, we may assume that $x$ can be expressed as such $K$-linear combinations as follows:
$$
\begin{array}{rcl}
x & = & k_{1}\alpha_1 + ... + k_{m}\alpha_m + k_{m+1}(\beta_{m+1}-\gamma_{m+1})+ ... + k_{m+p}(\beta_{m+p}-\gamma_{m+p})\\
&=& l_{1}\delta_1 + ... + l_{n}\delta_n + l_{n+1}(\zeta_{n+1}-\eta_{n+1}) + ... + l_{n+q}(\zeta_{n+q}- \eta_{n+q}),
\end{array}
$$
where $k_{1}, ..., k_{m+p}, l_{1}, ..., l_{n+q} \in K\backslash \{0\}$, $m, n, p, q$ are non-negative integers. Moreover, $\alpha_1$, ..., $\alpha_m$, $\beta_{m+1}-\gamma_{m+1}$, ...,  $\beta_{m+p}-\gamma_{m+p}$ are $K$-linearly independent special relations in $I_1$, and $\alpha_1$, ..., $\alpha_m$, $\beta_{m+1}$, ..., $\beta_{m+p}$, $\gamma_{m+1}$, ..., $\gamma_{m+p}$ are different paths. Similarly, $\delta_1$, ..., $\delta_n$, $\zeta_{n+1}-\eta_{n+1}$, ..., $\zeta_{n+q}-\eta_{n+q}$ are $K$-linearly independent special relations in $I_2$, and $\delta_1$, ..., $\delta_n$, $\zeta_{n+1}$, ..., $\zeta_{n+p}$, $\eta_{n+1}$, ..., $\eta_{n+p}$ are different paths.

For any $\alpha_i$, $1 \leq i \leq m$, we can easily have $$\alpha_i \in \{\delta_1, ..., \delta_n, \zeta_{n+1}, ..., \zeta_{n+q}, \eta_{n+1},...,\eta_{n+q}\}.$$
If $\alpha_i= \delta_j$ for some $j$ with $1\leq j \leq n$, then $k_i\alpha_i= l_j\delta_j \in I_1 \cap I_2$. If $\alpha_i= \zeta_{n+j}$ for some $j$ with $1\leq j \leq q$, then $\eta_{n+j} \neq \beta_{m+r}$ for any $1 \leq r \leq p$. Otherwise, the three different paths $\alpha_i= \zeta_{n+j}$, $\eta_{n+j}=\beta_{m+r}$, and $\gamma_{m+r}$ have the same source and the same target, which contradicts the condition that there are at most two paths from one vertex to another. Similarly, we also have $\eta_{n+j} \neq \gamma_{m+r}$ for any $1 \leq j \leq r$. Hence, we conclude $\eta_{n+j}= \alpha_s$ for some $s \neq i$ and $1 \leq s \leq m$. In this case, $k_i = l_{n+j}$ and $k_s= - l_{n+j}$. So we have $k_i \alpha_i+k_s\alpha_ s = l_j(\zeta_j - \eta_j)\in I_1 \cap I_2$.

Now we consider the special relation $\beta_{m+i} - \gamma_{m+i}$ for any $1\leq i \leq p$. It is clear that $$\beta_{m+i} \in \{\delta_1, ..., \delta_n, \zeta_{n+1}, ..., \zeta_{n+q}, \eta_{n+1},...,\eta_{n+q}\}.$$
If $\beta_{m+i}= \delta_j$ for some $1\leq j \leq n$, then $\gamma_{m+i}\neq \zeta_{n+r}$ for any $1\leq r \leq q$. Otherwise, the three different paths $\beta_{m+i}= \delta_j$, $\gamma_{m+i} = \zeta_{n+r}$, and $\eta_{n+r}$ have the same source and the same target, which is a contradiction. Similarly, we also have $\gamma_{m+i}\neq \eta_{n+r}$ for any $1 \leq r \leq q$. Thus, $\gamma_{m+i}= \delta_s$ for some $s \neq i$ and $1 \leq s \leq n$. In this case $k_{m+i} = l_j = - l_s$, and we have $k_{m+i}(\beta_{m+i}-\gamma_{m+i})=l_j \delta_j + l_s\delta_s \in I_1 \cap I_2$. If $\beta_{m+i} = \zeta_{n+j}$ for some $1 \leq j \leq q$, since there are at most two paths from one vertex to another, then clearly we have $\gamma_{m+i} = \eta_{n+j}$. Hence, $k_{m+i}= l_{n+j}$ and $k_{m+i}(\beta_{m+i}-\gamma_{m+i})=l_j(\zeta_{n+j}-\eta_{n+j})\in I_1 \cap I_2$.

Thus, we conclud that $x$ can be generated by special relations in $I_1 \cap I_2$. Thus $I_1 \wedge I_2 = I_1 \cap I_2$, which means that $\mathcal{SI}(KQ)$ is a sublattice of $\mathcal{I}(KQ)$. Note that $\mathcal{I}(KQ)$ is modular, so is $\mathcal{SI}(KQ)$.

$(2)\Rightarrow (3)\Rightarrow (4)$. It is clear that every modular lattice is strong lower semimodular and every strong lower semimodular lattice is lower semimodular.

$(4)\Rightarrow (1)$. Suppose that there are three different paths $\alpha$, $\beta$ and $\gamma$ between two vertices in the quiver $Q$. Then all the arrows and vertices that make up the three paths form a subquiver $Q'$. We claim that there is no arrow $a\in Q'$ such that $t(a)= s(\alpha)$. Suppose that there is an arrow $a\in Q'$ such that $t(a) = s(\alpha)$. Note that $a$ is a part of $\alpha$, or $\beta$, or $\gamma$. If $\alpha = u a v$ for some paths $u, v$ in $Q$, then it follows $s(ua)= s(\alpha)=t(a)=t(ua)$. Since $Q$ is acyclic, we obtain that $u = a = s(\alpha)$, that is, $a$ is a vertex, which is a contradiction. If $a$ is a part of $\beta$ or $\gamma$, the proof is similar. So we conclud that there is no arrow $a\in Q'$ such that $t(a)= s(\alpha)$. Dually, we also obtain that there is no arrow $a\in Q'$ such that $s(a)= t(\alpha)$.

Denote by $T$ the $K$-vector space with a basis the set of all paths in $Q$ but not in the subquiver $Q'$. It is easy to see that $T$ is a special algebraic ideal. Indeed, for any path $\xi$ in $T$ and $x, y$ in $S_Q$, if $x\xi y\neq 0$, then $x\xi y$ is a path not in the subquiver $Q'$, that is, $x\xi y \in T$.

Note that $span~\{\alpha, \beta-\gamma\} \cap T =\{0\}$. Let $I_1= span~\{\alpha, \beta-\gamma\} \oplus T$, the direct sum of the two $K$-vector spaces $span~\{\alpha, \beta-\gamma\}$ and $T$. We claim that $I_1$ is a special algebraic ideal. Indeed, for any paths $x, y$ in $Q$, if $x\alpha y= 0$ or $x\alpha y= \alpha$, then we have $x\alpha y \in I_1$. If $x\alpha y\neq 0$ and $x\alpha y\neq \alpha$, then $x\alpha y$ is a path longer than $\alpha$, we have shown that there is no arrow $a\in Q'$ such that $t(a) = s(\alpha)$, or $s(a)= r(\alpha)$. Thus, we have $x\alpha y\notin Q'$. That is, $x\alpha y \in T \subseteq I_1$.
Since $\beta$ and $\gamma$ have the same source and the same target, we know that $x\beta y= 0$ if and only if $x\gamma y=0$, $x\beta y= \beta$ if and only if $x\gamma y = \gamma$. If $x\beta y\neq 0$ and $x \beta y \neq \beta$, then $x\beta y$ is a path longer than $\beta$, so we have $x \beta y \in T$. Similarly, $x \gamma y \in T$. Thus we conclude $x(\beta- \gamma)y \in I_1$. This completed the proof that $I_1$ is a special algebraic ideal.

Similar to $I_1$, $I_2=span~\{\beta, \alpha-\gamma\}\oplus T$ is also a special algebraic ideal. One can see that $I_1 \vee I_2 = span~\{\alpha, \beta, \gamma\}\oplus T$, and $I_1 \cap I_2= span\{\alpha+ \beta -\gamma\}\oplus T$, which is not a special algebraic ideal. So it follows from the definition of the meet that $I_1 \wedge I_2 =T$. Thus, $I_1 \vee I_2 \succ I_1$ and $I_1 \vee I_2 \succ I_2$, but there exists a special algebraic ideal $I=span~\{\alpha\} \oplus T$ such that $I_1 \wedge I_2 \subsetneq I \subsetneq I_1$, which is a contradiction.
\end{prf}\medskip

\begin{eg} \label{eg3}
Let $Q$ be the Kronecker quiver in Figure \ref{Figure 6}.

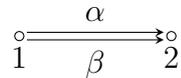
\begin{figure}[htbp]
\centering
\begin{tikzpicture}
\draw (1,1) circle(0.06);
\draw (3,1) circle(0.06);
\draw [-stealth] (1.1,1.05)--(2.9,1.05);\draw [-stealth] (1.1,0.95)--(2.9,0.95);
\node[below] at (1,1) {$1$};
\node[below] at (3,1) {$2$};
\node[above] at (2,1.05) {$\alpha$}; \node[below] at (2,0.95) {$\beta$};

\end{tikzpicture}
\caption{\footnotesize An example where $\mathscr{C}(S_Q)$ is modular.}\label{Figure 6}.
\end{figure}

The path semigroup $S_Q= \{0, e_1, e_2, \alpha, \beta\}$. There are 8 congruences: $\rho_1 = 1_{S_Q}$, $\rho_2$ is the congruence such that $(\alpha, \beta)\in \rho_2$ and other $\rho_2$-classes are trivial, $\rho_3$ is the Rees congruence generated by $\{0, \alpha\}$, $\rho_4$ is the Rees congruence generated by $\{0, \beta\}$, $\rho_5$ is the Rees congruence generated by $\{0, \alpha, \beta\}$, $\rho_6$ is the Rees congruence generated by $\{0, e_1, \alpha, \beta\}$, $\rho_7$ is the Rees congruence generated by $\{0, e_2, \alpha, \beta\}$, $\rho_8 = S_Q\times S_Q$.

See the Hasse diagram shown in Figure \ref{Figure 7}. The lattice $\mathscr{C}(S_Q)$ is modular but not distributive.

\begin{figure}[htbp]
\centering
\begin{tikzpicture}[scale=0.7]
\fill (0,0)circle(2pt); \node[below] at (0,0) {$1$};
\fill (0,1)circle(2pt); \node[right] at (0,1) {$2$};
\fill (0,2)circle(2pt); \node[above] at (0,2) {$5$};
\fill (0,4) circle(2pt); \node[above] at (0,4) {$8$};
\fill (-1,1) circle(2pt);   \node[left] at (-1,1) {$3$};
\fill (1,1) circle(2pt);   \node[right] at (1,1) {$4$};
\fill (-1,3) circle(2pt);   \node[left] at (-1,3) {$6$};
\fill (1,3) circle(2pt);   \node[right] at (1,3) {$7$};
\draw (0,0)--(-1,1);
\draw (0,0)--(0,1);
\draw (0,0)--(1,1);
\draw (-1,1)--(0,2);

\draw (0,1)--(0,2);
\draw (1,1)--(0,2);
\draw (0,2)--(1,3);
\draw (0,2)--(-1,3);
\draw (0,4)--(1,3);
\draw (0,4)--(-1,3);
\end{tikzpicture}
\caption{\footnotesize The Hasse diagram of the lattice $\mathscr{C}(S_Q)$ in Example \ref{eg3}.}\label{Figure 7}
\end{figure}

\end{eg}

 The following theorem provides equivalent conditions for the congruence lattice to be distributive.

\begin{thm}\label{distributive condition}
Let $Q$ be an acyclic quiver. Then the following statements are equivalent.
\begin{enumerate}[$(1)$]
\item There is at most one path from one vertex to another.
\item Each congruence $\rho \in \mathscr{C}(S_Q)$ is a Rees congruence.
\item $\mathscr{C}(S_Q)$ is distributive.
\end{enumerate}
\end{thm}

\begin{prf}
$(1)\Rightarrow (2)$. Let $\rho$ be any congruence. It is clear that $0\rho$ is a semigroup ideal of $S$. For any $(\alpha, \beta)\in \rho$, if $\alpha \in 0\rho$, then $\beta \in 0 \rho$, that is, $(\alpha, \beta) \in 0\rho \times 0\rho$. If $(\alpha, 0)\notin \rho$, then it follows from Lemma \ref{property of congruence} that $s(\alpha) = s(\beta)$ and $ t(\alpha) = t(\beta)$. Since $Q$ is acyclic and there is at most one path between any two vertices, we have $\alpha = \beta$. Therefore, we proved that $\rho = (0\rho \times 0\rho)\cup 1_{S_Q}$, which is a Rees congruence.

$(2)\Rightarrow (3)$. Since each congruence is a Rees congruence, the lattice $\mathscr{C}(S_Q)$ is isomorphic to the lattice of semigroup ideals. For any semigroup ideals $I_1$, $I_2$ in $S_Q$, it is clear that $I_1 \wedge I_2 = I_1 \cap I_2$ and $I_1 \vee I_2 = I_1 \cup I_2$. Thus $\mathscr{C}(S_Q)$ is distributive.

$(3)\Rightarrow (1)$. Suppose that there are two different paths $\alpha$ and $\beta$ between two vertices in the quiver $Q$. Then all arrows and vertices that make up the two paths form a subquiver $Q'$. We denote by $T$ the $K$-vector space whose basis is the set of all paths in $Q$ but not in the subquiver $Q'$. Then $T$ is a special algebraic ideal.

Note that the $K$-vector space $span~\{\alpha\} \cap T =\{0\}$. Let $I_1= span~\{\alpha\} \oplus T$, the direct sum of the two $K$-vector spaces. Similar to the proof of $(4)\Rightarrow (1)$ in Theorem \ref{modular condition}, $I_1$ is a special algebraic ideal of $KQ$.

Similarly, let $I_2= span~\{\alpha\} \oplus T$, and it is also a special algebraic ideal. Moreover, since $\alpha $ and $\beta$ have the same source and the same target, we know that $I_3= span~\{\alpha -\beta\} \oplus T$ is a special algebraic ideal as well. We can easily get that
$$ I_1 \wedge I_2 = I_1 \wedge I_3 = I_2 \wedge I_3 = T,$$ and $$I_1 \vee I_2 = I_1 \vee I_3 = I_2 \vee I_3 = span~\{\alpha,\beta\} \oplus T.$$
By Lemma \ref{diagrams of distri}, $\mathcal{SI}(KQ)$ is not distributive, which is a contradiction.
\end{prf}\medskip

For the path semigroups of acyclic quivers, it is proved that the modularity and distributivity of the congruence lattices are related to the number of paths between two vertices. Naturally, we have the following.
\begin{cor} Let $Q$ be an acyclic quiver and $\{Q_i\}_{i\in I}$ be the collection of all connected components of $Q$, where $I$ is an index set. Then we have the following.
\begin{enumerate}[$(1)$]
\item The lattice $\mathscr{C}(S_Q)$ is modular if and only if each lattice $\mathscr{C}(S_{Q_i})$ is modular for any $i\in I$.
\item The lattice $\mathscr{C}(S_Q)$ is distributive if and only if each lattice $\mathscr{C}(S_{Q_i})$ is distributive for any $i\in I$.
\end{enumerate}
\end{cor}

It follows from Gabriel's theorem that a connected, finite dimensional hereditary algebra is representation-finite if and only if the underlying graph of its quiver is one of the Dynkin graphs: $\mathbb{A}_m$ with $m\geq 1$; $\mathbb{D}_n$ with $n\geq 4$; and $\mathbb{E}_6$, $\mathbb{E}_7$, $\mathbb{E}_8$ (also in \cite[Theorem 5.10 of Chapter \uppercase\expandafter{\romannumeral7}]{Assem}).

In the end, we give the following corollary as an interesting application.

\begin{cor}\label{application tree}
Let $Q$ be a quiver. If the underlying graph of $Q$ is a tree, in particular, where it is a Dynkin graph, then $\mathscr{C}(S_Q)$ is distributive.
\end{cor}

\noindent {\bf Acknowledgements:} The authors would like to thank Dr. Baptiste Rognerud (Institut de Math\'{e}matiques de Jussieu - Paris Rive Gauche) for helpful suggestions and discussions.\\

\noindent {\bf Declaration:} The authors declare that there is no conflict of interest.

\end{document}